\documentclass[12pt]{amsart}
\input xypic
\usepackage{yfonts}

\usepackage{amssymb}
\newtheorem{thm}{Theorem}[section]
\newtheorem{lemma}[thm]{Lemma}

\newtheorem{prop}[thm]{Proposition}

\theoremstyle{definition}

\def\End{\operatorname{End}\nolimits}
\def\Ext{\operatorname{Ext}\nolimits}
\def\Hom{\operatorname{Hom}\nolimits}
\def\id{\operatorname{id}\nolimits}
\def\Ind{\operatorname{Ind}\nolimits}

\def\Gsl{{\mathfrak{sl}}}
\def\mqcoh{\operatorname{\!-qcoh}\nolimits}
\def\mMOD{\operatorname{\!-Mod}\nolimits}
\def\mModgr{\operatorname{\!-modgr}\nolimits}
\def\Res{\operatorname{Res}\nolimits}
\def\SL{\operatorname{SL}\nolimits}

\def\Spec{\operatorname{Spec}\nolimits}

\def\CC{{\mathcal{C}}}
\def\CL{{\mathcal{L}}}

\def\CN{{\mathcal{N}}}

\def\CO{{\mathcal{O}}}
\def\CT{{\mathcal{T}}}

\def\BA{{\mathbf{A}}}
\def\BG{{\mathbf{G}}}
\def\BL{{\mathbf{L}}}

\def\BP{{\mathbf{P}}}
\def\BQ{{\mathbf{Q}}}

\def\BT{{\mathbf{T}}}
\def\BZ{{\mathbf{Z}}}
\def\GS{{\mathfrak{S}}}

\def\iso{\buildrel \sim\over\to}

\begin{document}
\author{Rapha\"el Rouquier}
\address{UCLA Mathematics Department, Los Angeles, CA 90095-1555, USA.}
\email{rouquier@math.ucla.edu}
\thanks{The author gratefully acknowledges support from the NSF
(grant DMS-1702305) and from the Simons Foundation (grant \#376202).}

\title{$2$-Representations of $\Gsl_2$ from Quasi-maps}
\date{\today}
\maketitle

\section{Introduction}

There are two main classes of constructions of $2$-representations of
Kac-Moody algebras \cite{Rou2}. One is algebraic, for example via 
representations of cyclotomic quiver Hecke algebras, and the other uses
constructible sheaves, for example on quiver varieties. We describe here
a new type of $2$-representations, using coherent sheaves.

One of our motivations is to develop an affinization of the theory of
$2$-representations of Kac-Moody algebras.
One would want a theory of $2$-representations of
affinizations of symmetrizable Kac-Moody algebras, or rather of
the larger Maulik-Okounkov algebras \cite{MauOk}.

\smallskip
A classical theme of (geometric) representation theory is that affinizations
arise from degenerations.
The author and Shan, Varagnolo and Vasserot \cite{ShaVarVas,VarVas}
have proposed that the 
affinization of the monoidal category associated to the positive part of
a symmetric Kac-Moody algebra should be a full monoidal subcategory of the 
derived category of $\CO$-modules on the derived cotangent stack of
the moduli stack of representations of a corresponding quiver. A
description of this category by generators and relations is missing, even in
the case of $\Gsl_2$.

\smallskip

This article stems from efforts to better understand $2$-representations
on categories of coherent sheaves. A number of constructions have been
given by Cautis, Kamnitzer and Licata (cf e.g. \cite{CauKaLi}).
We study here a different
geometrical framework.
Feigin, Finkelberg, Kuznetsov, Mirkovi\'c and Braverman \cite{FeiFiKuMi,Bra}
have provided
a construction of Verma modules for complex semi-simple Lie algebras
using based quasi-map spaces from $\BP^1$ to flag varieties (zastavas).
We consider here the case of $\Gsl_2$, where the zastavas are smooth, and
are mere affine spaces. We show that coherent sheaves on zastavas provide a
$2$-Verma module for $\Gsl_2$ in the sense of Naisse-Vaz \cite{NaiVa1}.
Adding a superpotential and considering matrix factorizations, we obtain
a realization of simple $2$-representations of $\Gsl_2$.

\smallskip
I thank Sergei Gukov for useful discussions.

\section{Zastavas and correspondences}
\subsection{Quasi-maps}
We fix a field $k$ and we consider varieties over $k$.
The space of maps $\BP^1\to\BP^1$ of degree $d$
sending $\infty$ to $\infty$ identifies with the space
of pairs $(g,h)$ of polynomials such that
$g$ and $h$ have no common roots, $\deg(g(z)-z^d)<d$ and $\deg h<d$.

The zastava space $V_d$ of quasi-maps $\BP^1\to\BP^1$
defined in a neighborhood of $\infty$ and sending $\infty$ to $\infty$
is the space of pairs as above, without the condition on roots.
There is an isomorphism
$$\BA^{2d}\iso V_d,\
(a,b)\mapsto (g(z)=a_1+a_2z+\cdots+a_dz^{d-1}+z^d,
h(z)=b_1+b_2z+\cdots+b_dz^{d-1}).$$

There is an action of
$\BT=\BG_m\times \BG_m$ on $V_d$. The first $\BG_m$-action is by rescaling
the variable $z$ with weight $-2$. The second $\BG_m$-action
is by scalar action on $h(z)$ with weight $2$.

\subsection{Correspondences}
\label{se:coherent}
Let  $Y_d=V_d\times \BA^1$. We extend the $\BT$-action on $V_d$ to an
action on $Y_d$ by letting $\BT$ act on $\BA^1$ by weight $(-2,0)$.

We have a diagram of affine varieties with $\BT$-actions
$$\xymatrix{
	&	Y_d \ar[dl]_{(g,h,z_0)\mapsto (g,h)}^{\phi_d}
	\ar[dr]^{\ \ (g,h,z_0)\mapsto ((z-z_0)g,(z-z_0)h)}_{\psi_d} \\
	V_d &&
	V_{d+1}
	}$$

\smallskip
		We have functors
	$$F_d=\psi_{d*}\circ \phi_d^*:D^b_\BT(V_d\mqcoh)\to D^b_\BT(V_{d+1}\mqcoh)$$
	$$E_d=\phi_{d*}\circ \BL\psi_d^*:D^b_\BT(V_{d+1}\mqcoh)\to D^b_\BT(V_d\mqcoh).$$

	\subsection{Universal Verma module}
Let $U_v(\Gsl_2)$ be the quantum enveloping algebra of $\Gsl_2$. It is
the $\BQ(v)$-algebra generated by $e$, $f$ and $k^{\pm 1}$ subject to the 
relations
$$ke=v^2ek,\ kf=v^{-2}fk,\ ef-fe=\frac{k-k^{-1}}{v-v^{-1}}.$$

The universal Verma module $M_\kappa$ is the $U_v(\Gsl_2)$-module over
$\BQ(v,\kappa)$ with basis $(m_d)_{d\ge 0}$, with
$$k(m_d)=\kappa v^{-2d}m_d,\ e(m_d)=\delta_{d,0}m_{d-1} \text{ and }
f(m_d)=\frac{v^{d+1}-v^{-d-1}}{v-v^{-1}}\cdot\frac{\kappa v^{-d}-\kappa^{-1}v^{d}}{v-v^{-1}}$$

Specializing $\kappa$ to $v^{\lambda}$ gives the Verma module with
highest weight $\lambda$.

\subsection{Geometric realization}
\label{se:geometric}
Let $\CC$ be the category of bigraded vector spaces $N$
such that $\dim(\bigoplus_{i\le n,j\in\BZ}N_{ij})<\infty$
for all $n$. 
Taking bigraded
dimension gives an isomorphism 
$$K_0(\CC)\iso \BZ((v))\otimes \BZ[t^{\pm 1}],\
N\mapsto \sum_{i,j}v^it^j\dim N_{ij}.$$

We denote by $\CT_d$ 
the full triangulated subcategory of $D^b_{\BT}(V_d\mqcoh)$ 
generated by objects $N\otimes \CO_{V_d}$ for $N\in\CC$.

Given $N\in\CC$ with class $P\in \BZ((v))\otimes \BZ[t^{\pm 1}]$ and 
given $C\in D^b_{\BT}(V_d\mqcoh)$, we write $P\cdot C$ for the object
$N\otimes_k C$ of $D^b_{\BT}(V_d\mqcoh)$ (well defined up to isomorphism).

\medskip

We put $E=\bigoplus_{d\ge 0}E_d$ and
$F=\bigoplus_{d\ge 0}tv^{-2d}F_d[1]$. The following proposition
is an immediate consequence of Lemma \ref{le:EFonK0} below.
It is a variant of a result of Braverman and Finkelberg \cite{BraFi}.

\begin{prop}
	The actions of $[E]$ and $[F]$ on $M=\bigoplus_{d\ge 0}\BQ(v,t)
	\otimes_{\BZ[v,t]}
	K_0(\CT_d)$ give an action of $U_v(\Gsl_2)$ and there is
	an isomorphism of representations  
	$$\BQ((v))\otimes_{\BQ(v)}
	M_{tv^{-1}}\iso M,\ m_d\mapsto (1-v^2)^{d}[\CO_{V_d}].$$
\end{prop}

\subsection{Modules}
	Let $A_d=k[V_d]=k[a_1,\ldots,a_d,b_1,\ldots,b_d]$,
a bigraded algebra with $\deg(a_i)=(2(d-i+1),0)$ and $\deg(b_i)=(2(d-i+1),-2)$.

Let $B_d=k[Y_d]
	=k[a_1,\ldots,a_d,b_1,\ldots,b_d,c]$, a bigraded algebra
with $\deg(a_i)=(2(d-i+1),0)$, $\deg(b_i)=(2(d-i+1),-2)$ and $\deg(c)=(2,0)$.

	There is a bigraded  action of $A_d$ on $B_d$
	by multiplication and a bigraded action of $A_{d+1}$ on $B_d$ given by
	multiplication preceded by the morphism of algebras
	$$f:A_{d+1}\to B_d,\ 
a_i\mapsto a_{i-1}-ca_i \text{ and }
b_i\mapsto b_{i-1}-cb_i 
$$
where we put $a_0=b_0=b_{d+1}=0$ and $a_{d+1}=1$ in $B_d$.

\smallskip

	Via the equivalences $\Gamma:D^b_\BT(V_d\mqcoh)\iso 
	D^b_{bigr}(A_d\mMOD)$, the functors
	$E_d$ and $F_d$ become
	$$F_d=B_d\otimes_{A_d}-:D^b_{bigr}(A_d\mMOD)\to D^b_{bigr}(A_{d+1}\mMOD)$$
	$$E_d=B_d\otimes_{A_{d+1}}^{\BL} -:D^b_{bigr}(A_{d+1}\mMOD)\to
	D^b_{bigr}(A_d\mMOD).$$

\begin{lemma}
\label{le:EFonK0}
		We have $[E_d(A_{d+1})]=\frac{1}{1-v^2}[A_d]$ and
		$[F_d(A_d)]=\frac{(1-t^{-2}v^{2(d+1)})(1-v^{2(d+1)})}{1-v^2}[A_{d+1}]$.
\end{lemma}

\begin{proof}
	We have $B_d\simeq A_d\otimes k[c]$ as bigraded $A_d$-modules and
	the first statement follows.

	The second statement follows from Lemma \ref{le:Bd} below.
\end{proof}

\begin{lemma}
	\label{le:Bd}
	There is an exact sequence of bigraded $A_{d+1}$-modules
	$$0\to t^{-2}v^{2(d+1)}\frac{1-v^{2(d+1}}{1-v^2}A_{d+1}\to
	\frac{1-v^{2(d+1)}}{1-v^2}A_{d+1}\to B_d\to 0.$$
\end{lemma}

\begin{proof}
	Let $C=k[a_1,\ldots,a_d,c',b_1,\ldots,b_d,c'']$.
The morphism $f$ is the composition of the following morphisms of algebras:
	$$f_1:k[a_1,\ldots,a_{d+1},b_1,\ldots,b_{d+1}]\to
	k[a_1,\ldots,a_d,c',b_1,\ldots,b_{d+1}]$$
	$$b_i\mapsto b_i,\ a_i\mapsto 
\begin{cases}
	-c'a_1 & \text{ for }i=1\\
	a_{i-1}-c'a_i & \text{ for } 1<i\le d \\
	a_d-c' & \text{ for } i=d+1
\end{cases}$$
$$f_2:k[a_1,\ldots,a_d,c',b_1,\ldots,b_{d+1}]\iso
k[a_1,\ldots,a_d,c',b_1,\ldots,b_{d+1}]$$
$$a_i\mapsto a_i,\ c'\mapsto c',\ b_i\mapsto
	\begin{cases}
		b_i & \text { for } i\le d \\
		b_{d+1}+c' & \text{ for }i=d+1
	\end{cases}$$
$$f_3:k[a_1,\ldots,a_d,c',b_1,\ldots,b_{d+1}]\to C$$
$$a_i\mapsto a_i,\ c'\mapsto c',\
 b_i\mapsto 
\begin{cases}
	-c''b_1 & \text{ for }i=1\\
	b_{i-1}-c''b_i & \text{ for } 1<i\le d \\
	b_d-c'' & \text{ for } i=d+1
\end{cases}$$
$$f_4:C\to
k[a_1,\ldots,a_d,b_1,\ldots,b_d,c]$$
	$$a_i\mapsto a_i,\ b_i\mapsto b_i,\ c'\mapsto c,\ c''\mapsto c.$$

The first morphism makes $k[a_1,\ldots,a_d,c',b_1,\ldots,b_{d+1}]$ into
a free $k[a_1,\ldots,a_{d+1},b_1,\ldots,b_{d+1}]$-module with basis
	$(1,c',\ldots,c^{\prime d})$.

The third morphism makes $C$ into
a free $k[a_1,\ldots,a_d,c',b_1,\ldots,b_{d+1}]$-module with basis
	$(1,c'',\ldots,c^{\prime\prime d})$.

The last morphism makes $k[a_1,\ldots,a_d,b_1,\ldots,b_d,c]$ fit into
an exact sequence of
	$C$-modules
	$$0\to C \xrightarrow{c'-c''}C\xrightarrow{f_4}
k[a_1,\ldots,a_d,b_1,\ldots,b_d,c]\to 0.$$

	Let $L_1$ (resp. $L_0$)
	be the free $A_{d+1}$-module with basis $(e_{ij})_{0\le i,j\le d}$
	(resp. $(f_{ij})_{0\le i,j\le d}$).
	We have a commutative diagram of $A_{d+1}$-modules
	$$\xymatrix{
		0\ar[r]& L_1\ar[r]^-{d_1}\ar[d]_\sim^{\alpha_1} &
	L_0\ar[r]^-{d_0}\ar[d]_\sim^{\alpha_0} &
	k[a_1,\ldots,a_d,b_1,\ldots,b_d,c]\ar[r]\ar@{=}[d] & 0 \\
	0\ar[r]& C\ar[r]_{c'-c''}& C\ar[r]_-{f_4} &
	k[a_1,\ldots,a_d,b_1,\ldots,b_d,c]\ar[r] & 0 
	}$$
	where the structure of $A_{d+1}$-module on $C$ comes from 
	$f_3\circ f_2\circ f_1$, the one on $k[a_1,\ldots,a_d,b_1,\ldots,b_d,c]$
	from $f_4\circ f_3\circ f_2\circ f_1$ and where
	$$d_0(f_{ij})=c^{i+j}$$
	$$d_1(e_{ij})=
	\begin{cases}
		f_{i+1,j}-f_{i,j+1} & \text{ for }0\le i,j<d\\
		-(a_1f_{0j}+a_2f_{1j}+\cdots+a_{d+1}f_{dj}+f_{d,j+1}) &
		\text{ for }i=d,\ 0\le j<d \\
		b_1f_{i0}+b_2f_{i1}+\cdots+b_{d+1}f_{id} &
		\text{ for }j=d
	\end{cases}$$
	$$\alpha_1(e_{ij})=c^{\prime i}c^{\prime\prime j}$$
	$$\alpha_0(f_{ij})=c^{\prime i}c^{\prime\prime j}$$

	Note that $d_1$ and $d_0$ are morphisms of bigraded modules with
	$$\deg(e_{ij})=\begin{cases}
		(2(i+j+1),0) & \text{ for }j\neq d \\
		(2(i+d+1),-2) & \text{ for }j=d
	\end{cases}
	\text{ and }
	\deg(f_{ij})=(i+j,0).$$

	The $A_{d+1}$-module $L_0$ is generated by 
	$\{d_1(e_{ij})\}_{0\le i<d,0\le j\le d}$ and
	$(f_{i0})_{0\le i\le d}$. It follows that the complex
	$0\to L_1\xrightarrow{d_1}L_0\to 0$
	is homotopy equivalent to a complex of the form
$$0\to \bigoplus_{0\le j\le d}A_{d+1}e_{d,j}\to
\bigoplus_{0\le i\le d}A_{d+1}f_{i,d}\to 0.$$
The lemma follows.
\end{proof}

\section{$2$-Representations}
We construct now endomorphisms of $F$ and $F^2$, leading to a
structure of $2$-representation.

\smallskip
$\bullet\ $Let $\rho_d:Y_d\to\BA^1$ be the projection map. It provides
a morphim $k[X]=\Gamma(\CO_{\BA^1})\to \Gamma(\CO_{Y_d})$, hence a morphism
$k[X]\to\End(F_d)$.

\smallskip
$\bullet\ $There is an action of $\BG_a$ on $V_d\times V_{d+1}$ given by
$u\cdot \bigl((g,h),(g',h')\bigr)=\bigl((g,h),(g',h'+u)\bigr)$.
It provides a map $\mathrm{Lie}(\BG_a)=k\to
\Gamma(\CT_{V_d\times V_{d+1}})$ and
we denote by $\omega'$ the image of $1$, a vector field on $V_d\times V_{d+1}$.

The morphism $\phi_d\times \psi_d:Y_d\to V_d\times V_{d+1}$ is a closed
immersion and we identify $Y_d$ with its image. There is a canonical
isomorphism
\begin{equation}
	\label{eq:Faskernel}
	F_d\iso \pi_{2*}\bigl(\CO_{Y_d}\otimes
\pi_1^*(-)\bigr)
\end{equation}
where $\pi_1:V_d\times V_{d+1}\to V_d$ is the first projection and
$\pi_2:V_d\times V_{d+1}\to V_{d+1}$ the second projection.

Let $\omega''$ be the image of $\omega'$ by the composition of
canonical maps
$$\Gamma(\CT_{V_d\times V_{d+1}}\otimes\CO_{Y_d})\to
\Gamma(\CN_{Y_d/(V_d\times V_{d+1})})\iso \Ext^1_{\CO_{V_d\times V_{d+1}}}(
\CO_{Y_d},\CO_{Y_d}).$$
Via the isomorphism (\ref{eq:Faskernel}), $\omega''$ defines 
an element $\omega\in\Hom(F_d,F_d[1])$.

\smallskip
$\bullet\ $
There is an isomorphism
$$V_d\times\BA^2\iso Y_d\times_{V_{d+1}}Y_{d+1},\
((g,h,z_0,z'_0)\mapsto \bigl((g,h,z_0),((z-z_0)g,(z-z_0)h,z'_0)\bigr).$$

There is a commutative diagram
$$\xymatrix{
	&& V_d\times\BA^2 \ar[ddll]_{(g,h,z_0,z'_0)\mapsto (g,h)}
	\ar[d]_{\id\times\pi}\ar[ddrr]^{\ \ \ \ \ (g,h,z_0,z'_0)\mapsto
	((z-z_0)(z-z'_0)g,(z-z_0)(z-z'_0)h)} \\
	&& V_d\times(\BA^2/\GS_2) \ar[dll]^{\phi'_d}\ar[drr]_{\psi'_d} \\
	V_d&&&& V_{d+2}
}$$
for some maps $\phi'_d$ and $\psi'_d$, and where $\pi:\BA^2\to\BA^2/\GS_2$
is the quotient map.

Consequently, we obtain an isomorphism
\begin{equation}
\label{eq:FFaskernel}
F_{d+1}F_d\iso \psi'_{d*}((\CO_{V_d}\boxtimes\pi_*(\CO_{\BA^2}))
\otimes\phi'_d(-)).
\end{equation}

Let $\partial$ be the endomorphism of $k[X_1,X_2]=\Gamma(\CO_{\BA^2})$ given
by 
$$\partial(P)=\frac{P(X_1,X_2)-P(X_2,X_1)}{X_2-X_1}.$$
It induces an endomorphism of $\pi_*(\CO_{\BA^2})$, hence, via
the isomorphism (\ref{eq:FFaskernel}), an endomorphism $T$ of
$F_{d+1}F_d$.

\medskip
The data of $((E_d)_d,X,T)$ above gives rise to a $2$-representation of
$\Gsl_2^+$, but it does not extend to a $2$-representation of $\Gsl_2$.
But the data of $((E_d)_d,X,T,\omega)$ gives rise to a $2$-Verma module as
defined by Naisse and Vaz
(cf \cite{NaiVa1} and \cite[Definition 4.1]{NaiVa2}).

\begin{thm}
The functors $E_d$, $F_d$ and the natural transformations $X,\omega,T$
	define a $2$-Verma module for $\Gsl_2$
	on $\bigoplus_d\CT_d$ equivalent to the universal $2$-Verma module
	of \cite[\S 5.2]{NaiVa1}.
\end{thm}

\begin{proof}
	We show that our construction is equivalent to
the Naisse-Vaz universal Verma module \cite[\S 5.2]{NaiVa1}.

Let $\Omega_d=\Ext^*_{A_d}(A_d/(b_1,\ldots,b_d),A_d/(b_1,\ldots,b_d))$.
The canonical Koszul isomorphism of graded algebras
$\Lambda(b_1^*,\ldots,b_d^*)\iso \Ext^*_{k[b_1,\ldots,b_d]}(k,k)$
induces an isomorphism of graded algebras
$$\iota:k[a_1,\ldots,a_d]\otimes \Lambda(b_1^*,\ldots,b_d^*)\iso \Omega_d.$$
We denote by $\omega_i$ the image of $(-1)^{d+1-i}b_{d+1-i}^*$ in $\Omega_d$
and by $x_i$ the image of $(-1)^ia_{d+1-i}$. We have
$\Omega_d=k[x_1,\ldots,x_d]\otimes\Lambda(\omega_1,\ldots,\omega_d)$.

\medskip
Consider the morphism of algebras $h:A_d\otimes A_{d+1}\to B_d,\
r\otimes s\mapsto rf(s)$. We put $a_i=a_i\otimes 1$,
$b_i=b_i\otimes 1$, $a'_j=1\otimes a_j$ and $b'_j=1\otimes b_j$ for
$1\le i\le d$ and $1\le j\le d+1$.

The morphism $h$ is surjective, and its kernel is the ideal
generated by $\tilde{b}_i=b'_i+\tilde{c}b_i- b_{i-1}$
where $\tilde{c}=a_d-a'_{d+1}$.

Let $R=k[a_1,\ldots,a_d,\tilde{c}]\subset A_d\otimes A_{d+1}$ and 
$I=\bigoplus_{1\le i\le d+1}R\tilde{b}_i$, an $R$-submodule of
$L=\bigoplus_{1\le i\le d+1}Rb'_i\oplus
\bigoplus_{1\le i\le d}Rb_i$.
The morphism $h$ restricts to a surjective morphism of algebras
$S_R(L)\to B_d$ with kernel generated by $L$.

The orthogonal of $I$
in $\Hom_R(L,R)=\bigoplus_{1\le i\le d+1}Rb_i^{\prime *}\oplus
\bigoplus_{1\le i\le d}Rb_i^*$ is
$$I^\perp=\bigoplus_{1\le i\le d+1}R(b_{i+1}^{\prime *}+
b_i^*-\tilde{c}b_i^{\prime *}).$$

\smallskip
Define 
$$\Omega_{d,d+1}=\Ext^*_{A_d\otimes A_{d+1}}\Bigl(
\bigl(A_d/(b_1,\ldots,b_d)\bigr)\otimes
\bigl(A_{d+1}/(b_1,\ldots,b_{d+1})\bigr),B_d\Bigr).$$
The decomposition $A_d\otimes A_{d+1}=S_R(L)\otimes k[a'_1,\ldots,a'_d]$
induces an isomorphism
$$\Omega_{d,d+1}\iso 
\Ext^*_{S_R(L)}(
S_R(L)/(L),B_d)$$
hence 
$$\Omega_{d,d+1}\iso 
\Ext^*_{S_R(L)}( S_R(L)/(L),S_R(L)/(I))\iso \Lambda_R(\Hom_R(L,R)/I^\perp).$$
Let $\omega_i$ be the element of $\Omega_{d,d+1}$ corresponding to
the image of $(-1)^{d+2-i}b_{d+2-i}^{\prime *}$ in $\Hom_R(L,R)/I^\perp$.
Let $y_i$ be the image of $(-1)^ia_{d+1-i}$ in $\Omega_{d,d+1}$ and 
$\xi$ the image of $\tilde{c}$.
We have $\Omega_{d,d+1}=k[y_1,\ldots,y_d,\xi]\otimes\Lambda(\omega_1,\ldots,
\omega_{d+1})$.
The actions of $\Omega_d$ and $\Omega_{d+1}$ on $\Omega_{d,d+1}$ are given by 
multiplication preceded by morphisms of algebras
$$\Omega_d\to\Omega_{d,d+1},\ x_i\mapsto y_i,\ \omega_i\mapsto \omega_i+
\xi \omega_{i+1}$$
$$\Omega_{d+1}\to\Omega_{d,d+1},\ x_i\mapsto y_i+\xi y_{i-1},\ \omega_i\mapsto
\omega_i.$$

\smallskip
Let $\CT'_d$ be the full triangulated subcategory of $D^b_{bigr}(\Omega_d)$ 
generated by objects $N\otimes k[x_1,\ldots,x_d]$ for $N\in\CC$.

There is an equivalence of triangulated categories
$R\Hom^\bullet_{A_d}(A_d/(b_1,\ldots,b_d),-):\CT_d\iso \CT'_d$.
This equivalence intertwines the action of $E_d$ and $F_d$ with
the action of $\Omega_{d,d+1}\otimes_{\Omega_d}-$ and
$\Omega_{d,d+1}\otimes_{\Omega_{d+1}}-$. This shows our
construction is equivalent to that of Naisse and Vaz.
\end{proof}

\section{Finite-dimensional Simple modules}

We fix now $n\ge 0$. We define a simple $2$-representation of
$\Gsl_2$ by defining a superpotential on the universal $2$-Verma module 
and considering matrix factorizations.

\smallskip

Let $(g,h)\in V_d$. We have 
$$\frac{h(z^{-1})}{g(z^{-1})}=
z\frac{b_d+b_{d-1}z+\cdots+b_1z^{d-1}}{1+a_dz+\cdots+a_1z^d}=z
\sum_{i\ge 0} (b_d v_{i,d}+\cdots+b_1 v_{i-d+1,d})z^i$$
for some polynomials functions $v_{i,d}$ of $a_1,\ldots,a_d$ with
$v_{0,d}=1$ and $v_{i,d}=0$ for $i<0$.

\smallskip
We define a morphism $W_{d,n}:V_d\times\BA^n\to\BA^1$ by
$$W_{d,n}\bigl((g,h),(\gamma_1,\ldots,\gamma_n)\bigr)=
\sum_{i=0}^n \gamma_{i+1}(b_d v_{i,d}+\cdots+b_1 v_{i-d+1,d})$$
where we put $\gamma_{n+1}=1$.

\smallskip
Note that $W_{d+1,n}\circ(\psi_d\times \id)=W_{d,n}\circ(\phi_d\times\id)$
and we denote by $W_{d,n}$ that morphism $Y_d\times\BA^n\to\BA^1$.
We endow $\BA^n=\Spec k[\gamma_1,\ldots,\gamma_n]$ with an action of
$(\BG_m)^2$ with $\deg(\gamma_i)=(2(n+1-i),0)$.
This makes $W_{d,n}$ into a homogeneous map of degree $(2(n+1),-2)$.

\smallskip
We denote by $\CT_{d,n}$ the homotopy category of $(\BG_m)^2$-equivariant
matrix factorizations of $W_{d,n}$ on $V_d\times\BA^n$.
The functors $E_d$ and
$F_d$ of \S\ref{se:coherent} and \S\ref{se:geometric}
extend to functors between the categories
$\CT_{d,n}$ and $\CT_{d+1,n}$.

\begin{prop}
	We have $\CT_{d,n}=0$ if $d>n$.

	The action of $[E]$ and $[F]$ on $\bigoplus_{d=0}^{n-1}\BQ\otimes 
	K_0(\CT_{d,n})$ give an action of $U_q(\Gsl_2)$ and the corresponding
	representation is simple of dimension $n+1$.

	The data of $(E,F,T,X)$ define a $2$-representation of
	$\Gsl_2$ on $\bigoplus_{d=0}^n \CT_{d,n}$ equivalent to the homotopy
	category of bounded complexes of objects of the simple
	$2$-representation $\CL(n)$ of \cite[\S 4.3.2]{Rou2}.
\end{prop}

\begin{proof}
If $d>n$, then
$$W_{d,n}=
b_d(\gamma_1+\sum_{i=1}^{n} \gamma_{i+1}v_{i,d})+
b_{d-1}(\gamma_2+\sum_{i=2}^{n}\gamma_{i+1}v_{i-1,d})+\cdots+
b_{d-n+1}(\gamma_n+v_{1,d})+b_{d-n}.$$
As a consequence, the homotopy category of matrix factorizations
$\CT_{d,n}$ is $0$.

\smallskip
Assume now $d\le n$. We have
$$W_{d,n}=b_d(\gamma_1+\sum_{i=1}^{n} \gamma_{i+1}v_{i,d})+
b_{d-1}(\gamma_2+\sum_{i=2}^{n}\gamma_{i+1}v_{i-1,d})+\cdots+
b_1(\gamma_d+\sum_{i=d}^{n}\gamma_{i+1}v_{i-d+1,d}).$$

	Let $P_n=k[x_1,\ldots,x_n]$, a graded algebra with $\deg(x_i)=2$.
	Define
$$A_{d,n}=(A_d\otimes k[\gamma_1,\ldots,\gamma_n])/
	(\{b_r\}_{1\le r\le d}\bigcup\{\gamma_r+
	\sum_{i=r}^n\gamma_{i+1}v_{i-r+1,d}\}_{1\le r\le d}).$$
	The inclusion map induces an isomorphism
	$$k[a_1,\ldots,a_d]\otimes k[\gamma_{d+1},\ldots,
	\gamma_n]\iso A_{d,n}.$$
		Composing with the inverse of the isomorphism
	$$k[a_1,\ldots,a_d,\gamma_{d+1},\ldots,\gamma_n]\iso 
	P_n^{\GS_d\times\GS_{n-d}},\
	a_i\mapsto e_{d-i+1}(x_1,\ldots,x_d),\ \gamma_i\mapsto
	e_{n-i+1}(x_{d+1},\ldots,x_n)$$
we obtain an isomorphism of graded algebras
	$$P_n^{\GS_d\times\GS_{n-d}}\iso A_{d,n}.$$

	We view $A_{d,n}$ as a $\BZ$-graded algebra with $\deg(a_i)=2(d-i+1)$
	and $\deg\gamma_i=2(n+1-i)$.
	We have an equivalence $D^b(A_{d,n}\mModgr)\iso \CT_{d,n}$.

	Let $\CT'_{d,n}$ be the homotopy category of $(\BG_m)^2$-equivariant
	matrix factorizations of $W_{d,n}$ on $Y_d\times\BA^n$.
	Let $B'_{d,n}=A_{d,n}\otimes k[c]$. We have an equivalence
	$D^b(B'_{d,n}\mModgr)\iso \CT'_{d,n}$. 
	Let
	$B_{d,n}=
	B'_{d,n}\otimes_{A_{d+1}[\gamma_1,\ldots,\gamma_n]}A_{d+1,n}$.
	We have 
	$$B_{d,n}=B'_{d,n}/(\gamma_{d+1}+\sum_{i=d+1}^n
	\gamma_{i+1}\sum_{j=0}^{i-d}v_{j,d}c^{i-j-d})$$
	and the inclusion map induces an isomorphism
	$$k[a_1,\ldots,a_d,c,\gamma_{d+2},\ldots,\gamma_n]\iso B_{d,n}.$$
	Composing with the inverse of the isomorphism
	$$k[a_1,\ldots,a_d,c,\gamma_{d+2},\ldots,\gamma_n]\iso 
	P_n^{\GS_d\times\GS_{n-d-1}}$$
	$$a_i\mapsto e_{d-i+1}(x_1,\ldots,x_d),\ c\mapsto -x_{d+1},\
	\gamma_i\mapsto e_{n-i+1}(x_{d+2},\ldots,x_n)$$
we obtain an isomorphism of graded algebras
	$$P_n^{\GS_d\times\GS_{n-d-1}}\iso B_{d,n}.$$
There is a commutative diagram

$$\xymatrix{
	P_n^{\GS_{d+1}\times\GS_{n-d-1}} \ar@{^{(}->}[r] \ar[d]_\sim & 
	P_n^{\GS_d\times\GS_{n-d-1}} \ar[d]^\sim \\
	k[a_1,\ldots,a_{d+1},\gamma_{d+2},\ldots,\gamma_n]
	\ar[d]_\sim \ar[r] & k[a_1,\ldots,a_d,c,\gamma_{d+2},\ldots,
	\gamma_n] \ar[d]^\sim \\
	A_{d+1,n} \ar[r] & B_{d,n} \\
		A_{d+1}[\gamma_1,\ldots,\gamma_n] \ar[r]_{f}
		\ar@{->>}[u] &
		B_d[\gamma_1,\ldots,\gamma_n] \ar@{->>}[u]
}$$

	We deduce that there is a commutative diagram
	$$\xymatrix{
		&& \CT'_{d,n} \ar[ddll] \\
		&&D^b(B'_{d,n}\mModgr)\ar[u]^-\sim \ar[dddll]_-{\Res} \\
		 \CT_{d,n} 
		 && &&\CT_{d+1,n}\ar[uull] \\
		 && D^b(B_{d,n}\mModgr)\ar[uu]\ar[dll]^{\Res} \\
		D^b(A_{d,n}\mModgr)\ar[uu]^-\sim 
		&&&& D^b(A_{d+1,n}\mModgr)\ar[uu]_-\sim \ar[uuull]_-{B_{d,n}
		\otimes^\BL_{A_{d+1,n}}-} \ar[ull]^-{B_{d,n}
		\otimes^\BL_{A_{d+1,n}}-}\\
		&& D^b(P_n^{\GS_d\times\GS_{n-d-1}}\mModgr)\ar[uu]^\sim
		\ar[dll]^{\Res}
		\\
		D^b(P_n^{\GS_d\times\GS_{n-d}})\ar[uu]^-\sim &&&&
		D^b(P_n^{\GS_{d+1}\times\GS_{n-d-1}})\ar[uu]_-\sim
		\ar[ull]^-{\Ind}
	}$$

	It follows that $\CT_{d,n}$ is equivalent to the bounded
	derived category of finitely generated 
	graded $(k[a_1,\ldots,a_d]\otimes k[\gamma_{d+1},\ldots,
	\gamma_n])$-modules, where $\deg(a_i)=2(d-i+1)$ and 
	$\deg(\gamma_i)=2(n+1-i)$.
	Similarly, the homotopy category of $(\BG_m)^2$-equivariant
	matrix factorizations of $W_{d,n}$ on $Y_d\times\BA^n$ is
	equivalent to the bounded
	derived category of finitely generated 
	graded $(k[a_1,\ldots,a_d]\otimes k[\gamma_{d+1},\ldots,
	\gamma_n])$-modules, where $\deg(a_i)=2(d-i+1)$ and 
	$\deg(\gamma_i)=2i$. We recover the usual construction of
	the (homotopy category of the) simple $2$-representation $\CL(n)$ of
	$\Gsl_2$ (cf \cite[\S 5.2]{Rou1} and \cite[\S 4.3.2]{Rou2}).
\end{proof}

This construction is a Koszul dual counterpart of the construction of
\cite[\S 7]{NaiVa1} based on adding a differential.

\end{document}